\documentclass[12pt]{amsart}
\usepackage{latexsym, amsmath, amssymb, amsthm}
\usepackage{epsfig}
\usepackage{amscd}
\usepackage{latexsym}
\usepackage{pstricks,pst-node,cite}

\newtheorem{theorem}{Theorem}[section]

\newtheorem{corollary}[theorem]{Corollary}

\def\scaleddraw #1 by #2 (#3 scaled #4){{
  \dimen0=#1 \dimen1=#2
  \divide\dimen0 by 1000 \multiply\dimen0 by #4
  \divide\dimen1 by 1000 \multiply\dimen1 by #4
  \draw \dimen0 by \dimen1 (#3 scaled #4)}
  }

 \newcommand{\BC}{\mathbb{C}}

\newcommand{\Aut}{\mathrm{Aut}}
\newcommand{\C}[1]{\mathcal #1}

\ifx\blackandwhite\undefined
  \psset{linecolor=blue}\newrgbcolor{darkred}{.9 0 0}\newrgbcolor{emgreen}{0 .9 0}
  \newrgbcolor{lightgreen}{.7  .99 .7}\newrgbcolor{lightred}{.99 .7 .7} \else
  \psset{linecolor=blue}\newgray{darkred}{.7}\newgray{emgreen}{0.3}\newgray{pup}{.5}\newgray{xxx}{.5}
\fi
 \psset{linewidth=1pt,dash=3pt 3pt,doublesep=.05,dotsize=1pt 5}
\SpecialCoor

\begin{document}

\author{Dongseok Kim}
\address{Department of Mathematics \\ Kyungpook National University \\ Taegu, 702-201 Korea}
\email{dongseok@knu.ac.kr}
\thanks{}

\author{Hye Kyung Kim}
\address{Mathematics, Catholic University of Taegu\\Kyongsan,
712-702 Korea} \email{hkkim@cu.ac.kr}

\author{Jaeun Lee}
\address{Department of Mathematics \\Yeungnam University \\Kyongsan, 712-749, Korea}
\email{julee@yu.ac.kr}
\thanks{This research was supported by the Yeungnam University research grants in 2007.}

\subjclass[2000]{05C50, 05C25, 15A15, 15A18}
\keywords{generalized characteristic polynomials, graph bundles, the Bartholdi zeta
functions}

\title[Generalized characteristic polynomials of graph
bundles]{Generalized characteristic polynomials of graph bundles}

\begin{abstract}
In this paper, we find computational formulae for generalized characteristic polynomials of graph
bundles. We show that the number of spanning trees in a graph is
the partial derivative (at $(0,1)$) of the generalized characteristic
polynomial of the graph.
Since the reciprocal of the Bartholdi zeta function of a graph
can be derived from the generalized characteristic polynomial of a graph,
consequently, the Bartholdi zeta function of a graph bundle
can be computed by using our computational formulae.
\end{abstract}

\maketitle

\section{Introduction}

One of classical invariants in graph theory is the characteristic polynomial
which comes from the adjacency matrices. It displays not only graph theoretical
properties but also algebraic prospectives, such as spectra of graphs.
There have been many meaningful generalizations of characteristic polynomials~\cite{LVZ, CDS}.
In particular, we are interested in one found by Cvetkovic and alt. as a
polynomial on two variables ~\cite{CDS},
$$F_G(\lambda, \mu )=\det (\lambda I - (A(G) - \mu \mathcal{D}(G))).$$

The zeta functions of finite graphs \cite{Bar, Bass, Ih} feature of Riemann¡¯s zeta functions
and can be considered as an analogue of the Dedekind zeta
functions of a number field.
It can be expressed as the determinant of a perturbation of the Laplacian
and a counterpart of the Riemann hypothesis~\cite{ST}.
Bartholdi introduced the Bartholdi zeta function $Z_G(u,t)$ of a graph $G$ together with
a comprehensive overview and problems on the Bartholdi zeta functions~\cite{Bar}.
He also showed that the reciprocal of the Bartholdi zeta function of $G$ is
$$Z_G(u,t)^{-1}=\left(1-(1-u)^2t^2\right)^{\varepsilon_G-\nu_G}
\det\left[I-A(G)t+(1-u)(\mathcal{ D}_G-(1-u)I)t^2\right].$$
Kwak and alt. studied the Bartholdi zeta functions of
some graph bundles having regular fibers \cite{KLS}. Mizuno
and Sato also studied the zeta function and the Bartholdi zeta
function of graph coverings \cite{MS, MS2}. Recently, it was shown that
the Bartholdi zeta function $Z_G(u,t)$
can be found as the reciprocal of the generalized characteristic polynomials $F_G(\lambda, \mu)$
with a suitable substitution~\cite{KL07}.

The aim of the present article is to find computational
formulae for generalized characteristic polynomials $F_G(\lambda, \mu)$
of graph bundles and its applications. For computational formulae,
we show that if the fiber of the graph bundle is a Schreier graph,
the conjugate class of the adjacency matrix has a representative
whose characteristic polynomial
can be computed efficiently using the representation theory of
the symmetric group. We also provide
computational formulae for generalized characteristic polynomials of
graph bundles $G\times^\phi F$ where the images of $\phi$ lie in an
abelian subgroup $\Gamma$ of $\Aut(F)$.
To demonstrate the efficiency of our computation formulae,
we calculate the generalized characteristic polynomials $F_G$ of
some $K_n$-bungles $G\times^\phi K_n$. Consequently,
we can obtain the generalized characteristic polynomials $F_{K_{1,m}\times K_n}(\lambda, \mu)$ of
$K_{1,m}\times K_n$ which is a standard model of network with hubs. Its adjacency matrix, known as a ``kite",
is one of important examples in matrix analysis.

The outline of this paper is as follows. First, we review the
terminology of the generalized characteristic polynomials
and show that the number of spanning trees in a graph is the
partial derivative (at $(0,1)$) of the generalized characteristic
polynomial of the given graph in section~\ref{review}.
Next, we study a similarity of the adjacency  matrices of graph bundles
and find computational formulae for generalized characteristic polynomial
$F_G(\lambda, \mu)$ of graph bundles in section~\ref{bundle}.
In section~\ref{appl}, we find
the generalized characteristic polynomial of $K_{1,m}\times K_n$ and
find the number of spanning trees of $K_{1,m}\times K_n$.

\section{Generalized characteristic polynomials and complexity}\label{review}

In the section, we review the definitions and useful properties of
the generalized characteristic polynomials and find
the number of spanning trees in a graph using the generalized characteristic polynomials.

Let $G$ be an undirected finite simple graph with vertex set $V(G)$
and edge set $E(G)$. Let $\nu_G$ and $\varepsilon_G$ denote the
number of vertices and edges of $G$, respectively. An
\emph{adjacency matrix} $A(G)=(a_{ij})$ is the $\nu_G\times \nu_G$
matrix with $a_{ij}=1$ if $v_i$ and $v_j$ are adjacent and
$a_{ij}=0$ otherwise. The \emph{degree matrix}  $\mathcal{D}(G)$ of $G$ is
the diagonal matrix whose $(i,i)$-th entry is the degree $d_i^G=\mathrm{deg}_G(v_i)$ of
$v_i$ in $G$ for each $1\le i\le \nu_G$. The \emph{complexity} $\kappa(G)$
of $G$ is the number of spanning trees in $G$. An
\emph{automorphism} of $G$ is a permutation of the vertex set $V(G)$
that preserves the adjacency. By $|X|$, we denote the cardinality of a
finite set $X$. The set of  automorphisms forms a permutation group,
called the \emph{automorphism group} $\Aut(G)$ of  $G$. The
\emph{characteristic polynomial} of $G$, denoted by
$\Phi(G;\lambda)$, is the characteristic polynomial $\det(\lambda
I-A(G))$ of $A(G)$. Cvetkovic and alt. introduced a
polynomial on two variables of $G$,  $F_G(\lambda, \mu )=\det (
\lambda I - (A(G) - \mu \mathcal{D}(G)))$ as a generalization of
characteristic polynomials of $G$~\cite{CDS}, for example, the
characteristic polynomial of $G$ is $F_G(\lambda, 0 )$ and the
characteristic polynomial of the \emph{Laplacian matrix}
$\mathcal{D}(G) - A(G)$ of $G$ is $(-1)^{\nu_G}F_G(-\lambda, 1)$.

In~\cite{KL07}], it was shown that the Bartholdi zeta function
$Z_G(u,t)$ of a graph can be obtained from the polynomial $F_G(\lambda, \mu)$
with a suitable substitution as follows.
$$Z_G(u,t)^{-1}=\left(1-(1-u)^2t^2\right)^{\varepsilon_G-\nu_G}t^{\nu_G}F_G \left(\frac{1}{t}-(1-u)^2t, (1-u)t \right)$$
and
$$F_G(\lambda, \mu)=\frac{\lambda^{\nu_G}}{(1-\mu^2)^{\varepsilon_G}}Z_G
\left(1-\frac{\lambda\mu}{1-\mu^2},\frac{1-\mu^2}{\lambda}\right)^{-1}.$$

The complexities for various graphs have been studied \cite{MS, No}.
In particular, Northshield showed that the complexity of a graph $G$ can be
given by the derivative $$f'_G(1) = 2(\varepsilon_G -
\nu_G)~\kappa(G)$$ of the function $f_G(u) = \mathrm{det} [I -
u~A(G) + u^2~(\mathcal{D}(G)- I)]$, for a connected graph
$G$~\cite{No}. By considering the idea of taking the derivative, we find that the complexity of
a finite graph can be expressed as the partial derivative of
the generalized characteristic polynomial $F_G(\lambda, \mu)$ evaluated at $(0,1)$.

\begin{theorem} Let $F_G(\lambda, \mu )=\det ( \lambda I - (A(G) - \mu \mathcal{
D}(G)))$ be the generalized characteristic polynomial of a graph
$G$. Then the number of spanning trees in $G$, $\kappa(G)$, is
$$ \frac{1}{2 \varepsilon_G}~ \frac{\partial F_G}{\partial \mu} |_{(0, 1)},$$
where $\varepsilon_G$ is the number of edges of $G$. \label{spanning}\end{theorem}

\begin{proof} Let $B_{\lambda, \mu}= \lambda I - (A(G) - \mu \mathcal{
D}(G)) = (( b_{\lambda, \mu})_{ij})$ and let $B_{\lambda, \mu}^k =
(( b_{\lambda, \mu})^k_{ij})$ denote the matrix $B_{\lambda, \mu}$
with each entry of $k$-th row replace the corresponding partial
derivative with respect to $\mu$. Then
$$\begin{array}{ll} \displaystyle \frac{\partial}{\partial \mu}( \rm{det} ~B_{\lambda,
\mu}) =  \sum_{\sigma} \rm{sgn} (\sigma)~\frac{\partial}{\partial
\mu}~(\prod_i~((b_{\lambda, \mu})_{i \sigma(i)}) \\[3ex]
\displaystyle~~~~~~~~~~~~= \sum_k~\sum_{\sigma}
~\rm{sgn}~\prod_i~(b_{\lambda, \mu})^k_{i \sigma(i)}) = \sum_k~(
\rm{det} ~B_{\lambda, \mu}^k). \end{array}$$
 Since ~$(b_{0,
1})^k_{ij} = d^G_i~\delta_{ij} + a_{ij}(\delta_{kj}-1) = d_i
\delta_{ij} - a_{ij} + a_{ij}\delta_{ki},$
$$\frac{\partial F_G}{\partial \mu} |_{(0, 1)} = \prod_k ~\mathrm{det}~
(B_{0,1} + a_{ij}(\delta_{ki})).$$ Let
$(M^k)_{ij}=a_{ij}\delta_{ki}$ and let $(C_{0,1})_{ij} $ be the
cofactor of $b_{ij}$ in $B_{0,1}$. Then
$$\mathrm{det}~(B_{0,1} + M^k) = \sum_j (B_{0,1} + M^k)_{kj}(-1)^{k+j}
(C_{0,1})_{kj} = d_k~(C_{0,1})_{kk}.$$ Since
$\mathrm{det}~B_{0,1}=0 $, $(C_{0,1})_{ij}=\kappa(G)$ for all $i$
and $j$ \cite{Bi}. Hence,
$$\frac{\partial F_G}{\partial \mu} |_{(0, 1)} =
\sum_k~d_k~(C_{0,1})_{kk} = \kappa(G) \sum_k d_k = 2~\varepsilon_G
\kappa(G).$$ Thus,
$$\kappa(G) = \frac{1}{2 \varepsilon_G}~
\frac{\partial F_G}{\partial \mu} |_{(0, 1)}.$$
\end{proof}

\section{Generalized characteristic polynomials of graph bundles} \label{bundle}

Let $G$ be a connected graph and let $\vec{G}$ be the digraph
obtained from $G$ by replacing each edge of $G$ with a pair of
oppositely directed edges. The set of directed edges of $\vec{G}$ is
denoted by $E(\vec{G})$. By $e^{-1}$, we mean the reverse edge to an
edge $e\in E(\vec{G})$. We denote the directed edge $e$ of $\vec{G}$
by $uv$ if the initial and the terminal vertices of $e$ are $u$ and $v$,
respectively. For a finite group $\Gamma$, a $\Gamma$-\emph{voltage
assignment} of $G$ is a function $\phi:E(\vec{G})\rightarrow \Gamma$
such that $\phi(e^{-1})=\phi(e)^{-1}$ for all $e\in E(\vec{G})$. We
denote the set of all $\Gamma$-voltage assignments of $G$ by
$C^1(G;\Gamma)$.

Let $F$ be another graph and let $\phi\in C^1(G;\Aut(F))$. Now, we
construct a graph $G\times^\phi F$ with the vertex set
$V(G\times^\phi F)=V(G)\times V(F)$, and two vertices $(u_1,v_1)$
and $(u_2,v_2)$ are adjacent in $G\times^\phi F$ if either
$u_1u_2\in E(\vec{G})$ and $v_2=v_1^{\phi(u_1u_2)}=v_1\phi(u_1u_2)$
or $u_1=u_2$ and $v_1v_2\in E(F)$. We call $G\times^\phi F$ the
$F$-\emph{bundle over $G$ associated with $\phi$} (or, simply a
\emph{graph bundle}) and the first coordinate projection induces the
bundle projection $p^\phi:G\times^\phi F\rightarrow G$. The graphs
$G$ and $F$ are called the \emph{base} and the \emph{fibre} of the
graph bundle $G\times^\phi F$, respectively. Note that the map
$p^\phi$ maps vertices to vertices, but the image of an edge can be
either an edge or a vertex. If $F=\overline{K_n}$, the complement of
the complete graph $K_n$ of $n$ vertices, then an $F$-bundle over
$G$ is just an $n$-fold graph covering over $G$. If $\phi(e)$ is the
identity of $\Aut(F)$ for all $e\in E(\vec{G})$, then $G\times^\phi
F$ is just the Cartesian product of $G$ and $F$, a detail can be found in \cite{KL}.

Let $\phi$ be an $\Aut(F)$-voltage assignment of $G$. For each
$\gamma\in\Aut(F)$, let $\vec{G}_{(\phi,\gamma)}$ denote the
spanning subgraph of the digraph $\vec{G}$ whose directed edge set
is $\phi^{-1}(\gamma)$. Thus the digraph $\vec{G}$ is the
edge-disjoint union of spanning subgraphs
$\vec{G}_{(\phi,\gamma)}$, $\gamma\in\Aut(F)$. Let
$V(G)=\{u_1,u_2,\ldots,u_{\nu_G}\}$ and
$V(F)=\{v_1,v_2,\ldots,v_{\nu_F}\}$. We define an order relation
$\le$ on $V(G\times^\phi F)$ as follows: for
$(u_i,v_k),(u_j,v_\ell)\in V(G\times^\phi F)$, $(u_i,v_k)\le
(u_j,v_\ell)$ if and only if either $k<\ell$ or $k=\ell$ and $i\le
j$. Let $P(\gamma)$ denote the ${\nu_F}\times {\nu_F}$ permutation
matrix associated with $\gamma\in\Aut(F)$ corresponding to the
action of $\Aut(F)$ on $V(F)$, {\em i.e.}, its $(i,j)$-entry
$P(\gamma)_{ij}=1$ if $\gamma(v_i)=v_j$ and $P(\gamma)_{ij}=0$
otherwise. Then for any $\gamma,\delta\in\Aut(F)$,
$P(\delta\gamma)=P(\delta)P(\gamma)$.  Kwak and Lee expressed the
adjacency matrix $A(G\times^\phi F)$ of a graph bundle
$G\times^\phi F$ as follows.

\begin{theorem} [\cite{KL}] Let $G$ and $F$ be graphs and let $\phi$ be an
$\Aut(F)$-voltage assignment of $G$. Then the adjacent matrix of
the $F$-bundle $G\times^\phi F$ is
$$A(G\times^\phi
F)=\left(\sum_{\gamma\in\mathrm{ Aut}(F)}P(\gamma) \otimes
A(\vec{G}_{(\phi,\gamma)})\right) +A(F)  \otimes I_{\nu_G},$$
where $P(\gamma)$ is the ${\nu_F}\times {\nu_F}$ permutation
matrix associated with $\gamma\in\Aut(F)$ corresponding to the
action of $\Aut(F)$ on $V(F)$, and $I_{\nu_G}$ is the identity
matrix of order $\nu_G$.
\end{theorem}

For any finite group $\Gamma$, a \emph{representation} $\rho$ of a
group $\Gamma$ over the complex numbers is a group
homomorphism from $\Gamma$ to the general linear group $\mathrm{
GL}(r,\mathbb{C})$ of invertible $r\times r$ matrices over $\mathbb{
C}$. The number $r$ is called the {\em degree} of the representation
$\rho$~\cite{Sa}. Suppose that $\Gamma\le S_n$ is a permutation
group on $\Omega$. It is clear that $P: \Gamma\to
\mathrm{GL}(r,\mathbb{C})$ defined by $\gamma\to P(\gamma)$, where
$P(\gamma)$ is the permutation matrix associated with
$\gamma\in\Gamma$ corresponding to the action of~$\Gamma$ on
$\Omega$, is a representation of~$\Gamma$. It is called the
\emph{permutation representation}. Let
$\rho_1=1,\rho_2,\dots,\rho_\ell$ be the irreducible representations
of~$\Gamma$ and let $f_i$ be the degree of $\rho_i$ for each $1\le
i\le \ell$, where $f_1=1$ and $\sum_{i=1}^\ell f_i^2=|\Gamma|$. It
is well-known that the permutation representation $P$ can be
decomposed as the direct sum of irreducible representations :
$\displaystyle \rho =\oplus_{i=1}^\ell m_i \rho_i$ \cite{Sa}. In
other words, there exists an unitary matrix $M$ of order $|\Gamma|$
such that
$$M^{-1}P(\gamma)M=\bigoplus_{i=1}^\ell (I_{m_i} \otimes
\rho_i(\gamma))\eqno(2)$$ for any $\gamma\in\Gamma$, where $m_i\ge
0$ is the multiplicity of the irreducible representation $\rho_i$
in the permutation representation $P$ and $\sum_{i=1}^\ell
m_if_i=\nu_ F$. Notice that $m_1 \geq 1$ because it represents the number
of orbits under the action of the group $\Gamma$.

It is not hard to show that

$$(M \otimes I_{\nu_G})^{-1}A(G\times^\phi F)(M \otimes I_{\nu_G} ) =  \left[\bigoplus_{i=1}^{\ell} {
\sum_{\gamma\in\Gamma}I_{m_i} \otimes \rho_i(\gamma) \otimes
A(\vec{G}_{(\phi,\gamma)})}\right]  + A(F)\otimes I_{\nu_G}.$$

 For any vertex
$(u_i,v_k)\in V(G\times^\phi F)$, its degree is $d_i^G+d_k^F$,
where $d_i^G=\deg_G(u_i)$ and $d_k^F=\deg_F(v_k)$. Then, by our
construction of $\mathcal{ D}(G\times^\phi F)$, the diagonal
matrix $\mathcal{ D}(G\times^\phi F)$ is equal to $I_{\nu_F}
\otimes \mathcal{ D}(G)+\mathcal{D}(F) \otimes I_{\nu_G}$ and thus
$$(M \otimes I_{\nu_G})^{-1}\mathcal{ D}(G\times^\phi F)(
M\otimes I_{\nu_G}) = I_{\nu_F}\otimes \mathcal{ D}(G) + \mathcal{
D}(F)\otimes I_{\nu_G}.$$ Therefore, the matrix $A(G\times^\phi
F)-\mu \mathcal{ D}(G\times^\phi F)$ is similar to
$$\displaystyle{  \bigoplus_{i=1}^{\ell} I_{m_i} \otimes \left(\sum_{\gamma\in\Gamma}\rho_i(\gamma)
\otimes A(\vec{G}_{(\phi,\gamma)}) - I_{f_i} \otimes \mu \mathcal{
D}(G)\right)  } + \left(A(F) - \mu \mathcal{ D}(F) \right)\otimes
I_{\nu_G}.$$

By summarizing these, we obtain the following theorem.

\begin{theorem}\label{general}
Let $G$ and $F$ be two connected graphs and let $\phi$ be an
$\Aut(F)$-voltage assignment of $G$. Let $\Gamma$ be the  subgroup
of the symmetric group $S_n$. Furthermore, let
$\rho_1=1,\rho_2,\ldots,\rho_\ell$ be the irreducible
representations of $\Gamma$ having degree $f_1,f_2,\ldots, f_\ell$,
respectively. Then  the matrix $A(G\times^\phi F)-\mu \mathcal{
D}(G\times^\phi F)$ is similar to
$$\displaystyle{  \bigoplus_{i=1}^{\ell} I_{m_i} \otimes \left( \sum_{\gamma\in\Gamma} \rho_i(\gamma) \otimes A(\vec{G}_{(\phi,\gamma)})
  - I_{f_i} \otimes \mu \mathcal{ D}(G) \right)
} +  \left(A(F) - \mu \mathcal{ D}(F) \right) \otimes I_{\nu_G},$$
where $m_i\ge 0$ is the multiplicity of the irreducible
representation $\rho_i$ in the permutation representation $P$ and
$\sum_{i=1}^\ell m_if_i=\nu_ F$.\hfill\qed
\end{theorem}

A graph $F$ is called a \emph{Schreier graph} if there exists  a
subset $S$ of $S_{\nu_F}$ such that $S^{-1}=S$ and the adjacency
matrix $A(F)$ of $F$ is $\sum_{s\in S}P(s)$. We call such an $S$ the
\emph{connecting set} of the  Schreier graph $F$. Notice that a
Schreier graph with connecting set $S$ is a regular graph of degree
$|S|$ and most regular graphs are Schreier graphs \cite[Section
2.3]{GT}. The definition of the Schreier graph here is different
that of original one. But, they are basically identical\cite[Section
2.4]{GT}. Clearly, every Cayley graph is a Schreier graph.

\begin{theorem}\label{ge-bun} Let $G$ be a connected graph  and let $F$ be a Schreier graph
with connecting set $S$. Let $\phi:E(\vec{G})\to \mathrm{Aut}(F)$ be
a permutation voltage assignment. Let $\Gamma$ be the  subgroup of
the symmetric group $S_{\nu_F}$ generated by $\{\phi(e), s\,:\, e\in
E(\vec{G}), s\in S\}$. Furthermore, let
$\rho_1=1,\rho_2,\ldots,\rho_\ell$ be the irreducible
representations of $\Gamma$ having degree $f_1,f_2,\ldots, f_\ell$,
respectively. Then the matrix $A(G\times^\phi F)-\mu \mathcal{
D}(G\times^\phi F)$ is similar to
$$\displaystyle \bigoplus_{i=1}^\ell I_{m_i}\otimes
\left (\sum_{\gamma\in \mathrm{ Aut}(F)} \rho_i(\gamma) \otimes
A(\vec{G}_{(\phi,\gamma)})  - I_{f_i}\otimes \mu (\mathcal{
D}(G)+|S|I_{\nu_G})+ \left(\sum_{s\in S}\rho_i(s)\right)\otimes
I_{\nu_G}\right) ,$$ where $m_i\ge 0$ is the multiplicity of the
irreducible representation $\rho_i$ in the permutation
representation $P$ and $\sum_{i=1}^\ell m_if_i=\nu_ F$. \hfill\qed
\end{theorem}
\begin{proof}
Let $F$ be  a Schreier graph with a connecting set $S$. Then the
adjacency matrix of $F$ is $A(F) = \sum_{s \in S}~P(s)$. Hence, for
any  voltage assignment $\phi:E(\vec{G})\to \Aut(F)$, one can see
that
$$A(G\times^\phi
F)=\displaystyle \left(\sum_{\gamma\in\mathrm{
Aut}(F)}P(\gamma)\otimes A(\vec{G}_{(\phi,\gamma)}) \right)+
\sum_{s\in S}P(s)\otimes I_{\nu_G}.$$
 Let
$\Gamma$ be the subgroup of $S_{\nu_F}$ generated by $\{\phi(e),
s\,:\, e\in E(\vec{G}), s\in S\}$. Since $F$ is a regular graph of
degree $|S|$, one can see that
$$\mathcal{ D}(G\times^\phi F) = I_{\nu_F}\otimes (\mathcal{ D}(G) + |S|I_{\nu_G}).
$$  Now, one can have that
the matrix $A(G\times^\phi F)-\mu \mathcal{ D}(G\times^\phi F)$ is
similar to
$$\displaystyle\bigoplus_{i=1}^\ell I_{m_i}\otimes
\left (\sum_{\gamma\in \mathrm{ Aut}(F)}\rho_i(\gamma)\otimes
A(\vec{G}_{(\phi,\gamma)})  - I_{f_i}\otimes \mu (\mathcal{
D}(G)+|S|I_{\nu_G}) + \left(\sum_{s\in S}\rho_i(s)\right)\otimes
I_{\nu_G}\right).$$
\end{proof}

It is easy to see that the following theorem follows immediately from
Theorem~\ref{ge-bun}.

\begin{theorem}\label{c-ge-bun} Let $G$ be a connected graph  and let $F$ be a Schreier graph
with connecting set $S$. Let $\phi:E(\vec{G})\to \mathrm{Aut}(F)$ be
a permutation voltage assignment. Let $\Gamma$ be the  subgroup of
the symmetric group $S_{\nu_F}$ generated by $\{\phi(e), s\,:\, e\in
E(\vec{G}), s\in S\}$. Furthermore, let
$\rho_1=1,\rho_2,\ldots,\rho_\ell$ be the irreducible
representations of $\Gamma$ having degree $f_1,f_2,\ldots, f_\ell$,
respectively. Then  the characteristic polynomial $F_{G\times^\phi
F}(\lambda, \mu)$ of
 a graph bundle $G\times^\phi F$ is
{\small $$ \prod_{i=1}^\ell \mathrm{det}
\left[ I_{f_i}\otimes \left[(\lambda + \mu |S|) I_{\nu_G} + \mu
\mathcal{ D}(G)\right]  - \sum_{\gamma\in \mathrm{
Aut}(F)}\rho_i(\gamma) \otimes  A(\vec{G}_{(\phi,\gamma)})-
\left(\sum_{s\in S}\rho_i(s)\right) \otimes I_{\nu_G} \right]^{m_i},$$} where
$m_i\ge 0$ is the multiplicity of the irreducible representation
$\rho_i$ in the permutation representation $P$ and
$\sum_{i=1}^\ell m_if_i=\nu_ F$. \hfill\qed
\end{theorem}

Let $F=\overline{K_n}$ be the trivial graph on $n$ vertices.  Then
any $\Aut(\overline{K_n})$-voltage assignment is just a
permutation voltage assignment  defined in \cite{GT}, and
$G\times^\phi\overline{K_n}=G^\phi$ is just an $n$-fold covering
graph of $G$. In this case, it may not be a regular covering.
Now, the following comes from Theorem~\ref{general}.

\begin{corollary}\label{cover}Let $G$ be a connected graph  and let
$F=\overline{K_n}$. The characteristic polynomial
$F_{G^\phi}(\lambda, \mu)$ of the connected covering $G^\phi$ of a
graph $G$ derived from a permutation voltage assignment
$\phi:E(\vec{G})\to S_n $ is
{\small $$
F_G(\lambda, \mu)\times \prod_{i=2}^{\ell} \mathrm{det}\left[ I_{f_i} \otimes  \left[(\lambda + r)I_{\nu_G}+
\mu \mathcal{ D}(G)\right] - \sum_{\gamma\in\Gamma} \rho_i(\gamma)
\otimes A(\vec{G}_{(\phi,\gamma)}) - \left(\sum_{s\in S}
\rho_i(s)\right)\otimes I_{\nu_G}\right]^{m_i},$$}
where $m_i\ge 0$ is the multiplicity of the irreducible
representation $\rho_i$ in the permutation representation $P$ and
$\sum_{i=1}^\ell m_if_i=n$.\hfill\qed
\end{corollary}

Next, we consider the characteristic polynomial depending
on two variable of graph bundles $G\times^\phi F$ where the images
of $\phi$ lie in an abelian subgroup $\Gamma$ of $\Aut(F)$ and the
fiber  $F$ is $r$-regular. In this case, for any
$\gamma_1,\gamma_2\in\Gamma$, the permutation matrices $P(\gamma_1)$
and $P(\gamma_2)$ are commutative and $\mathcal{ D}_F=rI_{\nu_F}$.

It is well-known (see \cite{Bi}) that every permutation matrix
$P(\gamma)$ commutes with the adjacency matrix $A(F)$ of $F$ for
all $\gamma\in\Aut(F)$. Since the matrices $P(\gamma)$,
$\gamma\in\Gamma$, and $A(F)$ are all diagonalizable and commute
with each other, they can be diagonalized simultaneously. $i. e.$,
there exists an invertible matrix $M_\Gamma$ such that
$M_\Gamma^{-1}P(\gamma)M_\Gamma$ and $M_\Gamma^{-1}A(F)M_\Gamma$
are diagonal matrices for all $\gamma\in\Gamma$. Let
$\lambda_{(\gamma,1)},\ldots,\lambda_{(\gamma,{\nu_F})}$ be the
eigenvalues of the permutation matrix $P(\gamma)$ and let
$\lambda_{(F,1)},\ldots,\lambda_{(F,{\nu_F})}$ be the eigenvalues
of the adjacency matrix $A(F)$. Then
$$M_\Gamma^{-1}P(\gamma)M_\Gamma=\left[\begin{array}{ccc}
\lambda_{(\gamma,1)}&&{\bf 0}\\&\ddots&\\{\bf
0}&&\lambda_{(\gamma,{\nu_F})}\end{array}
\right]\quad\mbox{and}\quad
M_\Gamma^{-1}A(F)M_\Gamma=\left[\begin{array}{ccc}
\lambda_{(F,1)}&&{\bf 0}\\&\ddots&\\{\bf
0}&&\lambda_{(F,\nu_F)}\end{array}\right].$$ Using these similarities, we find that
 $$\begin{array}{l}
 (M_\Gamma \otimes I_{\nu_G})^{-1}\left(\displaystyle\sum_{\gamma\in
\Gamma}P(\gamma)\otimes A(\vec{G}_{(\phi,\gamma)}) +A(F)\otimes
I_{\nu_G}
\right)(M_\Gamma \otimes I_{\nu_G})\\[3ex]
\hspace{7cm}=\displaystyle\bigoplus_{i=1}^{\nu_F}\left(\displaystyle\sum_{\gamma\in
\Gamma}\lambda_{(\gamma,i)}A(\vec{G}_{(\phi,\gamma)})+
\lambda_{(F,i)}I_{\nu_G}\right).\end{array}$$
Recall that
$$\begin{array}{ll}(M \otimes I_{\nu_G})^{-1}\mathcal{ D}(G\times^\phi F)(M \otimes I_{\nu_G}
)\\[3ex]
 =I_{\nu_F} \otimes \mathcal{ D}(G)+r
I_{\nu_F}\otimes I_{\nu_G}
 = \displaystyle \bigoplus^{\nu_F}_{i=1}~(\mathcal{ D}(G)+rI_{\nu_G}) .
\end{array}$$

By summarizing these facts, we find the following theorem.

\begin{theorem}\label{ge-abl-bun} Let $G$  be a connected graph
and let $F$ be a connected regular graph of degree $r$. If the
images of $\phi \in C^1(G; \mathrm{Aut}(F))$ lie in an abelian
subgroup of $\mathrm{Aut}(F)$, then the matrix $A(G\times^\phi
F)-\mu \mathcal{ D}(G\times^\phi F)$ is similar to
$$\displaystyle \bigoplus^{\nu_F}_{i=1}\left((\lambda_{(F,i)} - r \mu)I_{\nu_G}
+ \left(\sum_{\gamma\in
\Gamma}\lambda_{(\gamma,i)}A(\vec{G}_{(\phi,\gamma)})
 - \mu \mathcal{ D}(G)\right)\right).$$
\hfill\qed
\end{theorem}

Notice that  the Cartesian product $G\times F$ of two graphs $G$
and $F$ is a $F$-bundle over $G$ associated with the trivial
voltage assignment $\phi$, {\em i.e.}, $\phi(e)=1$ for all $e\in
E(\vec{G})$ and $A(G)=A(\vec{G})$.  The following corollary comes
from this observation.

\begin{corollary}\label{ge-cart}
For any connected graph $G$ and a connected $r$-regular graph $F$,
the matrix $A(G\times F)-\mu \mathcal{ D}(G\times F)$ of   the
cartesian product $G\times F$ is similar to

$$\displaystyle \bigoplus_{i=1}^{\nu_F} \left((
\lambda_{(F,i)}-r\mu)I_{\nu_G}\! + (A(G)   -\! \mu \mathcal{ D}(G))
\right).$$ In particular, if $G$ is a regular graph of degree $d_G$,
then the matrix $A(G\times F)-\mu \mathcal{ D}(G\times F)$ of the
cartesian product $G\times F$ is
$$
\displaystyle \bigoplus_{i=1}^{\nu_F}
~\bigoplus_{j=1}^{\nu_G}\left( \lambda_{(F,i)} +
\lambda_{(G,j)}-(r+d)\mu\right),$$
 where
$\lambda_{(G,j)}$ $(1 \le j \le {\nu_G})$ and $\lambda_{(F,i)}$
$(1 \le i \le {\nu_F})$ are the eigenvalues of  $G$ and $F$,
respectively. \hfill\qed
\end{corollary}

Now, we consider the case the images of $\phi \in C^1(G; \mathrm{ Aut}(F))$
lie in an abelian subgroup of $\mathrm{Aut}(F)$.
 A vertex-and-edge weighted digraph
is a pair $D_{\omega}=(D, \omega)$, where $D$ is a digraph and
$\omega : V(D) \bigcup E(D) \rightarrow \mathbb{C}$ is a function.
We call $\omega$ the \emph{vertex-and-edge weight function} on $D$.
Moreover, if $\omega(e^{-1})= \overline{\omega(e)}$, the complex
conjugate of $\omega(e)$, for each edge $e \in E(D)$, we say that
$\omega$ is {\em symmetric}.
  Given any vertex-and-edge weighted digraph $D_{\omega}$, the
 adjacency matrix $A(D_{\omega})=(a_{ij})$
of $D_{\omega}$ is a square matrix of order $\mid V(D) \mid $
defined by
$$ a_{ij} =   \sum_{e \in E( \{ v_i \}, \{ v_j \})} \omega(e),$$
and the degree matrix $\mathcal{ D}_{D_\omega}$ is the diagonal
matrix whose $(i,i)$-th entry is $\omega(v_i)$. We define
$$F_{D_\omega}(\lambda, \mu)=\mathrm{
det}\left(\lambda I - \left(A(D_\omega) - \mu \mathcal{
D}_{D_\omega}\right)\right).$$ For any $\Gamma$-voltage assignment
$\phi$ of $G$, let
$\omega_i(\phi) : V(\vec G) \bigcup E(\vec G) \rightarrow \mathbb{C}$ be the
function defined by $$\omega_i (\phi)(v) = \mathrm{deg}_G(v), \quad\omega_i
(\phi)(e) = \lambda_{(\phi(e),i)}$$ for $e \in E(\vec G)$ and $v \in V(\vec G)$ where $i=1,2,\ldots,\nu_F$. Using
Theorem~\ref{ge-abl-bun}, we have the following theorem.

\begin{theorem}\label{c-abl-bun} Let $G$  be a connected graph
and let $F$ be a connected regular graph of degree $r$. If the
images of $\phi \in C^1(G; \mathrm{Aut}(F))$ lie in an abelian
subgroup of $\mathrm{Aut}(F)$, then the characteristic polynomial
$F_{G\times^\phi F}(\lambda, \mu)$ of
 a graph bundle $G\times^\phi F$ is
$$\displaystyle\prod_{i=1}^{\nu_F} ~F_{\vec G_{\omega_i(\phi)}} (\lambda+r\mu-\lambda_{(F,i)}, \mu).
$$\hfill\qed
\end{theorem}

Now, the following corollary follows immediately from Corollary
\ref{ge-cart}.

\begin{corollary}\label{c-cart}
For any connected graph $G$ and a connected $r$-regular graph $F$,
the characteristic polynomial $F_{G\times F}(\lambda, \mu)$ of the
cartesian product $G\times F$ is
$$\displaystyle \prod_{i=1}^{\nu_F} ~F_G( \lambda+r\mu-\lambda_{(F,i)}, \mu).$$
In particular, if $G$ is a regular graph of degree $d_G$, then the
characteristic polynomial $F_{G\times F}(\lambda, \mu)$ of the
cartesian product $G\times F$ is
$$
 \prod_{i=1}^{\nu_F}
~\prod_{j=1}^{\nu_G}\left(\lambda +
 (r + d)\mu-\lambda_{(G,j)} - \lambda_{(F,i)}\right),$$
 where
$\lambda_{(G,j)}$ $(1 \le j \le {\nu_G})$ and $\lambda_{(F,i)}$
$(1 \le i \le {\nu_F})$ are the eigenvalues of  $G$ and $F$,
respectively. \hfill\qed
\end{corollary}

\section{Generalized characteristic polynomial
of $K_{1,m}\times K_n$}\label{appl}

In this section, we find the generalized characteristic polynomial
of $K_{1,m}\times K_n$ and find the number of spanning trees of $K_{1,m}\times K_n$.
As we mentioned in introduction, $K_{1,m}\times K_n$ is a typical model for networks with hubs thus, we will count its spanning trees. Since $K_n$ features many nice structures, we discuss the generalized characteristic polynomial of graph bundles with a fiber, Cayley graph.

Let $\C{A}$ be a finite group with identity $id_\C{A}$ and let $S$ be
a set of generators for $\C{A}$ with the properties that $S =S^{-1}$
and $id_\C{A}\not\in S$, where $S^{-1}=\{x^{-1}\,|\,
x\in \Omega\}$. The \emph{Cayley graph} $Cay(\C{A}, S)$ is a
simple graph whose vertex-set and edge-set are defined as follows:
$$V(Cay(\C{A}, S))=\C{A}~\mathrm{and}~ E(Cay(\C{A}, S))=\{\{g, h\}\,
|\, g^{-1}h \in S\}.$$

From now on, we assume that $\C{A}$ is an abelian group of order $n$.
Let $G$ be a graph and let $\phi:E(\vec{G})\to \C{A}$ be an $\C{A}$-voltage assignment.
Notice that the left action $\C{A}$ on the vertex set $\C{A}$ of
$Cay(\C{A}, S)$ gives a group homomorphism from $\C{A}$ to $\Aut(Cay(\C{A}, S))$.
Let $P$ be the permutation representation of $\C{A}$ corresponding to the action.
Then the map $\tilde\phi:E(\vec{G})\to \Aut(Cay(\C{A}, S))$ defined by $\tilde\phi(e)=P(\phi(e))$ for any $e\in E(\vec{G})$
is an $\Aut(Cay(\C{A}, S))$-voltage assignment. We also denote it $\phi$.
Notice that every irreducible representation of an abelian group is linear. For convenience, let
$\chi_1$ be the principal character of $\C{A}$ and $\chi_2, \ldots, \chi_{n}$ be the other $n-1$ irreducible characters of $\C{A}$.
Now, by Theorem~\ref{ge-bun}, we have  that the matrix $A(G\times^\phi Cay(\C{A}, S))-\mu
\mathcal{D}(G\times^\phi Cay(\C{A}, S))$ is similar to
$$\begin{array}{l}
\displaystyle \left(|S|(1- \mu))I_{\nu_G}
+ \left(A(\vec{G})
 - \mu \mathcal{D}(G)\right)\right)\\[2ex]
 \hspace{1cm} \displaystyle \oplus~~
 \bigoplus^{n}_{i=2}\left(\left(\chi_i(S) - |S| \mu\right)I_{\nu_G}
+ \left(\sum_{\gamma\in
\C{A}}\chi_i(\gamma)A(\vec{G}_{(\phi,\gamma)})
 - \mu \mathcal{D}(G)\right)\right),\end{array}$$
where $\chi_i(S)=\sum_{s\in S}\chi_i(s)$ for each $i=2,3,\ldots,n$.
By Theorem~\ref{c-abl-bun}, we have that the characteristic polynomial
$F_{G\times^\phi Cay(\C{A}, S)}(\lambda, \mu)$ of
 a graph bundle $G\times^\phi Cay(\C{A}, S)$ is
$$F_G (\lambda+ |S| (\mu-1), \mu) \times \prod_{i=2}^{n} ~F_{\vec G_{\omega_i(\phi)}} (\lambda+ |S| \mu- \chi_i(S), \mu),$$
where
$\omega_i(\phi) : V(\vec{G}) \bigcup E(\vec{G}) \to \BC$ be the
function defined by $$\omega_i (\phi)(v) = \mathrm{deg}_G(v), ~~~~\omega_i
(\phi)(e) = \chi_i(\phi(e))$$ for $e \in E(\vec G)$ and $v \in V(\vec G)$ where $i=2, 3, \ldots, n$.
Let $K_n$ be the complete graph on $n$ vertices. Then $K_n$ is isomorphic to $Cay(\C{A},\C{A}-\{id_\C{A}\})$
for any group $\C{A}$ of order $n$. Since $\chi_1(\C{A}-\{id_\C{A}\})=n-1$ and $\chi_i(\C{A}-\{id_\C{A}\})=-1$
for each $i=2,3,\ldots,n$, we have
$$F_{G\times^\phi K_n}(\lambda, \mu)=F_G (\lambda+ (n-1) (\mu-1), \mu) \times \prod_{i=2}^{n} ~F_{\vec G_{\omega_i(\phi)}} (\lambda+ (n-1) \mu +1, \mu).$$
Moreover over if $\phi$ is the trivial voltage assignment, then
$$ \displaystyle F_{G\times K_n}(\lambda, \mu)=F_G (\lambda+ (n-1) (\mu-1), \mu)\times F_G(\lambda+ (n-1) \mu +1, \mu)^{n-1}.$$

Let $G$ be the complete bipartite graph $K_{1,m}$ which also called  a star graph. Notice that $K_{1,m}$ is a tree
and hence every graph bundle $K_{1,m}\times ^\phi F$ is isomorphic to the cartesian product $K_{1,m}\times  F$
of $K_{1,m}$ and $F$. It is known \cite{KL07} that for any natural numbers $s$ and $t$
$$F_{K_{s,t}}(\lambda,\mu)=(\lambda+t\mu)^{s-1}(\lambda+s\mu)^{t-1}\left[(\lambda+s\mu)(\lambda+t\mu)-st\right]$$
and hence $F_{K_{1,m}}(\lambda,\mu)=(\lambda+\mu)^{m-1}\left[(\lambda+\mu)(\lambda+m\mu)-m\right]$.
Now, we can see that
$$\begin{array}{l}
\displaystyle F_{K_{1,m}\times^\phi K_n}(\lambda, \mu)\\[1ex]
\hspace{1cm}= \displaystyle F_{K_{1,m}} (\lambda+ (n-1) (\mu-1), \mu) \times  ~F_{K_{1,m}}(\lambda+ (n-1) \mu +1, \mu)^{n-1} \\[1ex]
\hspace{1cm} =  \displaystyle [\lambda+n\mu-(n-1)]^{m-1}\left\{[\lambda+n\mu-(n-1)][\lambda+(m+n-1)\mu-(n-1)]-m\right\} \\[1ex]
  \hspace{2cm} \times \displaystyle [\lambda+n\mu+1]^{(m-1)(n-1)}\left\{[\lambda+n\mu+1][\lambda+(m+n-1)\mu+1]-m\right\}^{n-1}. \end{array}
      $$

Now, by applying Theorem~\ref{spanning}, we have that the number of spanning trees of
$K_{1,m}\times K_n$ is

$$\kappa(K_{1,m}\times K_n) = n^{n-2}(m+n+1)^{n+1}(n+1)^{(m-1)(n-1)}.$$

\end{document}